\documentclass{amsart}

\issueinfo{00}
  {}
  {}
  {2004}

\copyrightinfo{2004}
  {}

\usepackage[T1]{fontenc}    
\usepackage{a4wide}
\usepackage{amsmath}        
\usepackage{amssymb}
\usepackage{amsthm}
\usepackage{paralist}
\usepackage{graphicx}
\usepackage{color}
\usepackage{bbm}            
\usepackage{url}            
\usepackage{calc}           
\usepackage[mathscr]{eucal}
\usepackage{ae,aecompl}    
\usepackage{xspace}        
\theoremstyle{plain}
\newtheorem{lemma}{Lemma}[section]
\newtheorem{theorem}[lemma]{Theorem}

\newtheorem{corollary}[lemma]{Corollary}
\newtheorem{proposition}[lemma]{Proposition}
\theoremstyle{definition}
\newtheorem{definition}[lemma]{Definition}

\theoremstyle{remark}

\newcommand{\R}{\mathbbm{R}}

\newcommand{\Z}{\mathbbm{Z}}

\newcommand{\eps}{\varepsilon}
\renewcommand{\phi}{\varphi}
\DeclareMathOperator{\aff}{aff}

\renewcommand{\vec}[1]{{\boldsymbol {#1}}}
\newcommand\vv{\vec{v}}
\newcommand\ww{\vec{u}}
\newcommand\uu{\vec{w}}
\newcommand\xx{\vec{x}}
\newcommand\yy{\vec{y}}
\newcommand\bb{\vec{b}}
\newcommand\cc{\vec{c}}
\newcommand\ee{\vec{e}}
\newcommand\nn{\vec{n}}
\newcommand\zero{\vec{0}}

\graphicspath{{pictures/}}

\begin{document}

\title{Projected Products of Polygons}%
\author{G\"unter M.~Ziegler}
\address{Inst.\ Mathematics, MA 6-2, TU Berlin, D-10623 Berlin, Germany}
\email{ziegler@math.tu-berlin.de}
\thanks{Partially supported by Deutsche
  Forschungs-Gemeinschaft (DFG), via the
  DFG Research Center ``Mathematics in the Key Technologies'' (FZT86),
  the Research Group ``Algorithms, Structure, Randomness'' (Project ZI 475/3),
  and a Leibniz grant (ZI 475/4)}
\date{\today}
 

\begin{abstract}\noindent
We construct a $2$-parameter family of $4$-dimensional polytopes
$\pi(P^{2r}_n)$ with extreme combinatorial structure: In this family, 
the ``fatness'' of the $f$-vector gets arbitrarily close to~$9$,
the ``complexity'' (given by the flag vector) gets arbitrarily close to~$16$.

The polytopes are obtained from suitable deformed products of even
polygons by a projection to four-space.
\end{abstract}

\maketitle

\section{Introduction}

According to Steinitz' paper \cite{Stei3} from 1906, the 
cone spanned by the $f$-vectors of
$3$-dimensional polytopes (with apex
at the $f$-vector of the $3$-simplex) is given by
\[
f_2\le 2f_0-4\qquad\textrm{and}\qquad f_0\le2f_2-4
\]
in conjunction with Euler's formula, $f_1=f_0+f_2-2$;
see also \cite[Sect.~10.3]{Gr1-2}.
Moreover, the cone of flag vectors of $3$-polytopes is 
also given by these two inequalities, 
together with the usual Dehn--Sommerville relations,
such as $f_{02}=2f_1$ and $f_{012}=4f_1$.
Here the extreme cases, polytopes whose
$f$- or flag vectors lie at the boundary of the cone,
are given by the 
simplicial polytopes (for which Steinitz' first inequality is tight) 
and the simple polytopes (second inequality tight).
Moreover, \emph{all} the integer points that satisfy the
above conditions do correspond to actual polytopes.

For $4$-dimensional polytopes, such a complete and simple
answer is not to be expected; 
compare Gr\"unbaum \cite[Sect.~10.4]{Gr1-2},
Bayer \cite{Bay} and H\"oppner \& Ziegler \cite{Z59}.
However, we may now be getting close
to a full description of the convex cone spanned by the $f$-vectors
(with apex at the $f$-vector of the $4$-simplex): In the projective coordinates
introduced in~\cite{Z82},
\[
\phi_0\ :=\ \frac{f_0-5}{f_1+f_2-20}\qquad\textrm{and}\qquad
\phi_3\ :=\ \frac{f_3-5}{f_1+f_3-20},
\]
we can write the known necessary conditions as 
\[
\phi_0\ge0,\quad\phi_3\ge0,\qquad
\phi_0+3\phi_3\le1,\quad3\phi_0+\phi_3\le1,\qquad\textrm{and}\qquad
\phi_0+\phi_3\le\tfrac25,
\]
where the first two conditions are trivial,
the second two have simplicial resp.\ simple polytopes as extreme cases,
and the last condition is a lower bound that follows from
``$g_2^{\mathrm{tor}}\ge0$.''
This is indeed a complete description of the cone
if and only if there are polytopes for which the sum
$\phi_0+\phi_3$ is arbitrarily small, that is, for which the
\emph{fatness} parameter
\[
F(P)\ :=\ \frac1{\phi_0+\phi_3}\ =\ \frac{f_1+f_3-20}{f_0+f_2-10}
\]
is arbitrarily large \cite{Z82}.
This observation has sparked a certain race for ``fat'' $4$-dimensional
polytopes. The following table summarizes the main steps.
Most of the examples that appear there are 2-simple and 2-simplicial,
with a symmetric $f$-vector; the first infinite family
of such polytopes was constructed by Eppstein, Kuperberg \& Ziegler
\cite{Z80}; a simple construction appears in Paffenholz \& Ziegler
\cite{Z89}. We call a polytope ``even'' if its graph is bipartite.
\[
\begin{tabular}{|l|l|l|l|l|}
\hline
4-polytopes & fatness      & property                & reference & date\\
\hline\hline
simple or simplicial            & $<\,3$     & &&\\ 
products    & $3-\eps$     & simple     &&\\ 
$24$-cell   & 4.526        & 2-simple, 2-simplicial  & 
                                            Schl\"afli~\cite{Schla} &1852\\
dipyramidal 720-cell& 5.020& 2-simple, 2-simplicial  & 
                                               Gev\'ay~\cite{Gevay} &1991\\
neighborly cubical 
            & $5-\eps$     & even    & Joswig \& Ziegler \cite{Z62} &2000\\
{[E-construction]} & 5.048 & 2-simple, 2-simplicial  & 
                                        Eppstein et al.\ \cite{Z80} &2003\\
$E(C_m\times C_n)$& $6-\eps$ & 2-simple, 2-simplicial  &
                                    Paffenholz \cite{paffenholz-pc} &2004\\
Projected products & $9-\eps$ & even       & here & \\
\hline
\end{tabular}
\]
A flag vector parameter
that is similar to fatness, called \emph{complexity},
was also introduced in~\cite{Z82}.
This is an invariant of the flag vector, defined by
\[
C(P)\ :=\ \frac{f_{03}}{f_0+f_3-10}\ =\ \frac{g_2(P)}{g_1(P)+g_1(P^*)}+3,
\]
where $g_1=f_0-5$ and $g_2=f_{03}-3f_0-3f_3+10$
are components of the toric $g$-vector of the polytope,
while $g_1(P^*)=f_3-5$ and $g_2(P)=g_2(P^*)$ refer to 
$g$-entries of the dual polytope; compare \cite{Sta7}
and~\cite{kalai87:_rigid_i}.
All $4$-polytopes satisfy $C(P)\ge3$.
Fatness and complexity are roughly within a factor of~$2$:
$C(P)\le 2F(P)-2$ and $F(P)\le 2C(P)-2$.
In particular, it is not known whether $C(P)$ can be arbitrarily
large.
Previously, the polytopes with the largest known 
complexity were the ``neighborly cubical polytopes''
of Joswig \& Ziegler \cite{Z62}, of complexity $8-\eps$.
Our present construction yields ``neighborly cubical polytopes''
in the special case $n=4$
(possibly special cases of those of \cite{Z62},
with a simpler description), but for $n,r\rightarrow\infty$
it yields complexity values as large as $16-\eps$.

The concept of ``strictly preserving a face'' used in the following
theorem will be explained in Section~\ref{sec:projections}.
(Compare the concept of faces in the ``shadow boundary'' of
a projection, e.g.\ in~\cite{Bar4}.)

\begin{theorem}\label{thm:ppp}
Let $n\ge4$ be even and $r\ge2$.
Then there is a $2r$-polytope $P_n^{2r}\subset\R^{2r}$, 
combina\-torially equivalent
to a product of $r$ $n$-gons, $P_n^{2r}\cong (C_n)^r$,
such that the projection $\pi:\R^{2r}\rightarrow\R^4$
to the last four coordinates
preserves the $1$-skeleton as well as 
all the ``$n$-gon $2$-faces'' of $P_n^{2r}$.
\end{theorem}

In the following, we describe in brief the main ingredients
for the construction and sketch the proof for its correctness.
Detailed proofs, the complete combinatorial characterization of the
resulting polytopes, possible extensions, further remarkable aspects
(such as the polyhedral surfaces of high genus embedded in the
$2$-skeleta of the resulting $4$-polytopes; cf.\ 
\cite{schroeder04}) as well as necessity
of the restrictions imposed here (e.g., that $n$ must be even)
are topics of current research and will be presented later.

\subsubsection*{Acknowledgements}
The intuition for the construction given here 
grew from previous joint work and current discussions
with Nina Amenta, Michael Joswig, and Thilo Schr\"oder.

\section{Products and Deformed Products}\label{sec:products}

The combinatorial structure of the
products of polygons $(C_n)^r$ is easy to describe:
These are simple $2r$-polytopes, with $f_0=n^r$ vertices,
$f_1=rn^r$ edges, and $f_{2r-1}=rn$ facets.
In general, its non-empty faces are products of non-empty faces
of the polygons, so 
\[
\sum_{i=0}^{2r-1}f_{2r-i}t^i\ \ =\ \ (1+nt+nt^2)^r.
\]
The $2$-dimensional faces of~$(C_n)^r$, and thus
of any polytope combinatorially equivalent to~$(C_n)^r$,
may be split into two classes: There are 
$rn^{r-1}$ faces that are $n$-gons,
which we refer to \emph{polygons}, which arise as products of one
of the $n$-gons with a vertex from each of the other factors;
and there are $\binom r2 n^r$ \emph{quadrilaterals} that 
(in $(C_n)^r$) arise
as products of edges from two of the factors with vertices from
the others. 
Thus, in total $(C_n)^r$ has $f_2=r n^{r-1} + \binom r2 n^r$ $2$-faces.

In the case $n=4$, the polygon $2$-faces of $(C_n)^r$
are $4$-gons, but we nevertheless treat the $r4^{r-1}$ polygons
and the $\binom{r}2 4^r$ quadrilaterals separately in this case.

An inequality description for such a product polytope may be 
obtained as 
\[
\left(\raisebox{-39mm}{\input{pictures/matrix0a.pstex_t}}\!\right)\xx 
\ \ \le\ \ 
\left(\raisebox{-39mm}{\input{pictures/vector0a.pstex_t}}\!\right),
\]
assuming that $V\xx\le\bb$ is a correct description
for an $n$-gon: For this it is necessary and sufficient
that the row vectors $\vv_i$ of $V$ are non-zero and distinct
and that they positively span~$\R^2$,
that the components $b_i$ of~$\bb$ are positive,
and that the rescaled vectors $\tfrac1{b_i}\vv_i$ are in convex
position (the vertices of the polar of the polygon).

For this we say that a finite set of
vectors $\vv_1,\dots,\vv_k\in\R^k$ \emph{positively spans} if 
it satisfies the following equivalent conditions:
\begin{compactenum}[(i)]
\item every vector $\xx\in\R^d$ is a linear combination of 
the vectors $\vv_i$, with non-negative coefficients,
\item every vector $\xx\in\R^d$ is a linear combination of 
the vectors $\vv_i$, with positive coefficients,
\item the vectors $\vv_i$ span $\R^d$, and
      $\zero\in\R^d$ is a linear combination of 
      the vectors $\vv_i$, with positive coefficients
      (that is, the vectors $\vv_i$ are \emph{positively dependent}).
\end{compactenum}
In the following, we will need ``deformed products''
(as described in Amenta \& Ziegler  \cite{Z51a}) of polygons.
For this, we look at systems of the form
\[
\left(\raisebox{-39mm}{\input{pictures/matrix1.pstex_t}}\!\right)\xx 
\ \ \le\ \ 
\left(\raisebox{-39mm}{\input{pictures/vector1a.pstex_t}}\!\right).
\]
Given any such left-hand side matrix for such a system,
we can adapt the right-hand side
so that the resulting polytope is \emph{combinatorially equivalent}
to~$(C_n)^r$: For this all components of $\bb_k$ have to be sufficiently
large compared to $\bb_1,\dots,\bb_{k-1}$,
for $k=2,3,\dots,r$. (Compare \cite{Z51a} and \cite{Z62}.)

\section{Projections}\label{sec:projections}

We work with a  rather restrictive concept of
faces ``being preserved under projection.''

\begin{definition}[Strictly preserving faces under projection]%
\label{def:strictly}%
Let $\pi:P\rightarrow Q=\pi(P)$ be a projection of polytopes.
Then a face $G\subseteq P$ is 
\emph{strictly preserved} if
\begin{compactenum}[\rm(i)]
\item its image $\pi(G)$ is a face of~$Q$,
\item the map $G\rightarrow\pi(G)$ is a bijection, and
\item the pre-image of the image is $G$, that is, $\pi^{-1}(\pi(G))=G$.
\end{compactenum}
\end{definition}

In the definition, Conditions (ii) and (iii) are both needed.
Indeed, in the projection ``to the second coordinate'' displayed in our figure, 
the vertex $v$ is strictly preserved, 
but the vertex~$w$ and the edge $e$ are not: 
For $w$ condition (iii) fails, while for $e$ condition (ii) is
violated. 
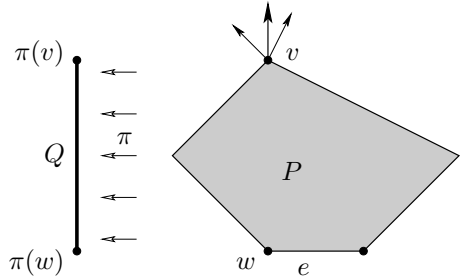
\[
\input{pictures/project.pstex_t}
\]

For simplicity, the following characterization result
is given only in a coordinatized version, for the
projection ``to the last $d$ coordinates.''

We say that a vector $\cc$ \emph{defines} the face 
$G\subseteq P$ given by all the points of $P$ that have
maximal scalar product with $\cc$. This describes exactly
all the vectors in the relative interior of the
\emph{normal cone} of~$G$. If $P$ is full-dimensional,
this interior of the normal cone consists of all the
positive combinations of outer facet normals $\nn=\nn_F$ to the facets
$F\subset P$ that contain~$G$.
(Compare \cite[Sects.~2.1, 3.2, 7.1]{Z35}.)

\begin{proposition}\label{prop:preserve}%
Let $\pi:\R^{e+d}\rightarrow\R^d$, $(\xx',\xx'')\mapsto\xx''$ be the
projection to the last~$d$ coordinates, and
let $P\subset\R^{e+d}$ be an $(e+d)$-dimensional polytope,
and let $G$ be a face of~$P$.
Then the following three conditions are equivalent:
\begin{enumerate}[\rm(1)]
\item
$G$ is strictly preserved by the projection $\pi:P\rightarrow\pi(P)=Q$.
\item
Any $\cc'\in\R^e$ arises as the first $e$ components 
of a vector $(\cc',\cc'')$ that defines~$G$.
\item
The vectors $\nn'$, given by the first $e$ components 
of the normal vectors $(\nn',\nn'')=\nn_F$ to facets $F$ of~$P$ that
contain $G$, positively span $\R^e$.
\end{enumerate}
\end{proposition}

\begin{proof}
Here we only establish ``$(3)\Rightarrow(1)$,'' which is used
in the following.

If the vectors $\nn'$ are positively dependent,
then some positive combination of the vectors $(\nn',\nn'')=\nn_F$ 
yields $(\zero,\cc'')=:\cc$.
A point $\xx\in P$ lies in the face $G\subseteq P$ if and only
if its scalar product with each facet normal $\nn$
is maximal. This happens if and only if
$\cc^t\xx$ is maximal, that is, iff
$(\cc'')^t\xx''$ is maximal under the restriction $\xx''\in\pi(P)$.
Thus we have established that under the assumption (3),
$\pi(G)=:\bar G$ is a face of~$\pi(P)$, and 
$\pi^{-1}(\pi(G))=\pi^{-1}(\bar G)=G$; that is,
Conditions (i) and (iii) of Definition~\ref{def:strictly}
are satisfied.

Since the vectors $\nn''$ additionally are positively spanning,
we know that if $\xx=(\xx',\xx'')\in G$, then
$(\xx',\xx''+\yy'')\notin G$ for $\yy''\neq\zero$:
Every such $\yy''$ arises as a positive combination
of the vectors $\nn''$, and since the scalar product with the corresponding
combination of the $(\nn',\nn'')$ is maximized over $P$
at $(\xx',\xx'')$, a scalar product with $(\xx',\xx''+\yy'')$
will be larger than the maximum over~$P$.

This implies that the projection
$\aff(G)\rightarrow\pi(\aff(G))=\aff(\pi(G))$
is injective. In particular $G\rightarrow\pi(G)$ is injective,
and this establishes part~(ii) of Definition~\ref{def:strictly}.
\end{proof}

\section{Construction}

\begin{proposition}[Construction for the proof of Theorem~\ref{thm:ppp}]
For $n\ge4$ even and $r\ge2$, let $P_n^{2r}$ be defined by
the linear inequality system
\[
\left(\raisebox{-55mm}{\input{pictures/matrix2.pstex_t}}\right)\xx 
\ \ \le\ \ 
\left(\raisebox{-55mm}{\input{pictures/vectorgrey.pstex_t}}\right).
\]
Here the left-hand side coefficient matrix 
$A_{n,r}^\eps\in\R^{rn\times 2r}$
contains blocks of size $n\times 2$, where
\[
V=\raisebox{-8mm}{\input{pictures/matrixVgrey.pstex_t}}
\longrightarrow
V^\eps=\raisebox{-8mm}{\input{pictures/matrixVepsgrey.pstex_t}},
\quad
W=\raisebox{-8mm}{\input{pictures/matrixWgrey.pstex_t}},
\quad
U=\raisebox{-8mm}{\input{pictures/matrixUgrey.pstex_t}}\ \in\ \R^{n\times2},
\]
with
\[
\vv_0=(1,0),\  
\vv_1=(0,0)=\zero,\  
\uu_0=(0,1),\  
\uu_1=(-3,-\tfrac23),\  
\ww_0=(-\tfrac{31}4,\tfrac12),\  
\ww_1=(9,-\tfrac23).
\]
The block $V^\eps$ arises from $V$ by an $\eps$-perturbation:
\[
\vv_i^\eps\ =\ \begin{cases}
\hphantom{\eps}\big(1-\eps (n-2-2i)^2,\eps (n-2-2i)\big)
&\textrm{for }i=0,2,4,\ldots,n-2,\\
         {\eps}\big(1-\eps (n-2-2i)^2,\eps (n-2-2i)\big)
&\textrm{for }i=1,3,5,\ldots,n-3,\\
\eps(-1,0)\ =\ (-\eps,0)                     
&\textrm{for }i=n-1
\end{cases}
\]
for a sufficiently small $\eps>0$.
All entries of $A^\eps_{n,r}$ outside the $r+(r-1)+(r-2)=3r-3$ blocks
of types $V^\eps$, $W$ and $U$ are zero.

Let the right-hand side vector be such that $\bb_1$ is given by 
$b_{1,i}=1$ for even $i$, and $b_{1,i}=\eps$ for odd~$i$,
and by $\bb_k=M^{k-1}\bb_1$ for sufficiently large $M$.

Then $P_n^{2r}$ has the properties claimed by Theorem~\ref{thm:ppp}.
In particular, it is a deformed product of $r$ $n$-gons, and
all its polygon $2$-faces survive the projection to 
the last $4$ coordinates.
\end{proposition}

\begin{proof}
The rows $\vv_i^\eps$ of $V^\eps$ are
indeed in cyclic order:
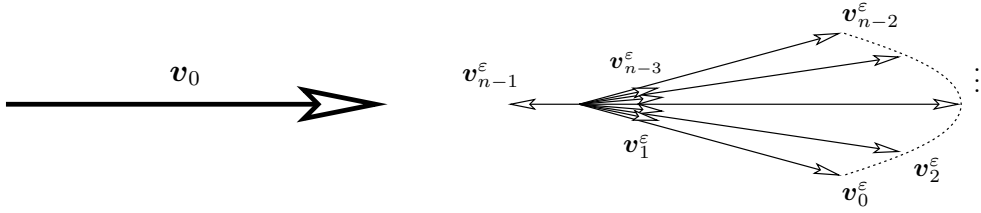
\[
\input{pictures/perturb.pstex_t}
\]
Moreover, rescaled as 
$\frac1{b_{k,i}}\vv_i^\eps=\frac1{M^{k-1}\eps}\vv_i^\eps$ for odd~$i$ and as
$\frac1{b_{k,i}}\vv_i^\eps=\frac1{M^{k-1}    }\vv_i^\eps$ for even $i$ 
they are in convex position, if $\eps$ is small; so
 $V^\eps\xx\le\bb_i$ defines a convex $n$-gon.
Thus for sufficiently small $\eps$ and sufficiently large~$M$,
the polytope $P_{2r}$ is indeed a deformed product of polygons,
as discussed in Section~\ref{sec:products}.

Now we show that for sufficiently small $\eps$,
all the polygon $2$-faces of $P_n^{2r}$ survive the projection
to the last $4$ coordinates. For this, we verify that the 
left-hand side matrix with $V$-blocks instead of $V^\eps$-blocks,
which we denote by $A_{n,r}^0=A_{n,r}$,
satisfies the linear algebra condition
dictated by Proposition~\ref{prop:preserve}(3).
This is sufficient, since the ``positively spanning'' condition is
stable under perturbation by a small~$\eps$.

Any polygon $2$-face $G$ of the simple $2r$-polytope $P^{2r}_n$ is defined
by  the facet normals to the $2r-2$ facets that contain~$G$.
The facet normals correspond to the rows of the inequality system,
and thus for the facet normals of a polygon $2$-face 
one has to choose two cyclically adjacent rows from each block
(corresponding to a vertex from each factor polygon), 
except from one of the blocks no row is taken.
Moreover, due to the structure of the matrices $U$, $V$, and $W$,
in which rows alternate, any choice of two cyclically-adjacent
rows from a block yields the same pair of rows (only 
the order is not clear, but it also does not matter).

Thus, to apply Proposition~\ref{prop:preserve}(3) we have to show:
\begin{quote}\emph{%
If one of the $r$ pairs of rows is deleted from the reduced matrix
\[
A'_{n,r}\ \ =\ \ 
\left(
\raisebox{-25mm}{\input{pictures/matrix_reduced1.pstex_t}}%
\right)
\ \ \in\ \R^{2r\times (2r-4)},
\]
then the remaining $2r-2$ rows\\
{\rm(a)} \ span $\R^{2r-4}$, and\\
{\rm(b)} \ have a linear dependence with strictly positive coefficients.}
\end{quote}
\smallskip

Let us establish (b) first. 
For this, let 
\[
\alpha_k\ :=\ 2^k+2^{-k}-2\qquad\textrm{and}\qquad
\beta_k \ :=\ 2^k+\tfrac54 2^{-k}-\tfrac94.
\]
These sequences are designed to be non-negative,
$\alpha_k,\beta_k\ge0$ for all $k\in\Z$, with equality only for $k=0$.
Thus for (b) it suffices to verify
\begin{quote}\emph{%
  For any $1\le t\le r$, the rows of $A'_{n,r}$ 
  are positively dependent with coefficients $\alpha_{k-t}$ for the
  even-index row from the $k$-th block, and $ \beta_{k-t}$ for the
  odd-index row from the $k$-th block.}
\end{quote}
since the (two) vectors in the $k$-th block thus get
zero coefficients, so they may be
deleted from any linear dependence (with otherwise positive
coefficients).
Thus we are led to the condition
\[
\alpha_{k-1}\vv_0 + 
\alpha_k    \uu_0 + \beta_k\uu_1 + 
\alpha_{k+1}\ww_0 + \beta_{k+1}\ww_1 \ =\ \zero,
\]
which is needed to hold for $k\le |r-2|$, but which we impose for all
$k\in\Z$.
The choice of vectors $\vv_0,\uu_0,\uu_1,\ww_0,\ww_1$ 
is designed to satisfy this condition. Indeed,
except for the choice of a basis, which we took to be  
$\vv_0=(1,0)$ and $\uu_0=(0,1)$, the 
configuration of five vectors $\vv_0,\uu_0,\uu_1,\ww_0,\ww_1$
are uniquely determined by the condition. 
\smallskip

For Property (a), we have to show that if one of the $r$ pairs of rows
is deleted from the matrix $A'_{n,r}$, 
then the resulting matrix still has full rank.
If the first or the second pair of rows is deleted,
then we still have the last $2r-4$ rows, and they
form a block upper triangular matrix, which has full
rank since its diagonal block
\[
\left(
\raisebox{-2.5mm}{\input{pictures/matrix_rankU.pstex_t}}%
\right)
\]
is non-singular.
If a later pair of rows is deleted, then we are faced
with the task to show that the $2k\times 2k$ matrices 
$M_k$ of the form
\[
M_k\ \ :=\ \ 
\left(
\raisebox{-15.5mm}{\input{pictures/matrix_reduced2a.pstex_t}}%
\right)
\ \ \in\ \R^{2k\times 2k}
\]
are non-singular. To verify this
(without proving explicitly that
$\det M_k=\frac{(2^k-1)^2}{3^k}$, which
might need combinatorial ingenuity)
we use our knowledge about row combinations
of $M_k$. Indeed, if we sum the rows of $M_k$ with
coefficients 
$(\alpha_0,\beta_0,\alpha_1,\ldots,\alpha_{k-1},\beta_{k-1})$,
then this will result in the linear combination
of the three rows of the matrix
\[
H_3\ \ =\ \ 
\left(
\raisebox{-4mm}{\input{pictures/matrix_rank3.pstex_t}}%
\right)\ \in\ \R^{3\times2r}
\]
with the coefficients $(-\alpha_{-1},-\alpha_k,-\beta_k)$,
since $\vv_1=\zero$. Similarly, 
if we sum the rows of $M_k$ with coefficients 
$(\alpha_1,\beta_1,\alpha_2,\ldots,\alpha_k,\beta_k)$,
then we get a linear combination of the same three rows,
with coefficients $(-\alpha_0,-\alpha_{k+1},-\beta_{k+1})$.
And if we use coefficients
$(\alpha_2,\beta_2,\alpha_3,\ldots,\alpha_{k+1},\beta_{k+1})$
to sum the rows of $M_k$, then the result will be a sum 
with coefficients $(-\alpha_1,-\alpha_{k+2},-\beta_{k+2})$.
The coefficient matrix
\[
\left(
\begin{smallmatrix}
-\alpha_{-1} & -\alpha_k & -\beta_k\\
-\alpha_0 & -\alpha_{k+1} & -\beta_{k+1}\\
-\alpha_1 & -\alpha_{k+2} & -\beta_{k+2}
\end{smallmatrix}
\right)
\]
is non-singular for $k\ge0$: Its determinant is
$\tfrac38(2^k-1+2^{-k-2})$.
Thus the full row-space
of $H_3$ is contained in the row space of $M_k$.
In particular, we find the unit vectors
$\ee_{2k-1},\ee_{2k}\in\R^{2k}$
in the row space of $H_3$, and thus of $M_k$, and
this allows us to complete the argument by induction.
\end{proof}

\section{Flag vectors}\label{sec:flag}

\begin{proposition}\label{cor:flag}
The $4$-polytope $\pi(P_n^{2r})$
has the flag vector 
\begin{eqnarray*}
(f_0,f_1,f_2,f_3;f_{03}) & = & 
(n^r,rn^r, \tfrac54rn^r-\tfrac34n^r+rn^{r-1}, 
           \tfrac14rn^r-\tfrac12n^r+rn^{r-1}; 4rn^r-4n^r)\\
&=&(4n, 4rn,   5rn-3n+4r, rn-2n+4r; 16rn-16n)\cdot \tfrac14 n^{r-1}
\end{eqnarray*}
\end{proposition}

\begin{proof}
We obtain $f_0=n^r$ and $f_1=rn^r$ from the
products $(C_n)^r$, which are simple $2r$-polytopes
with $n^r$ vertices. With the abbreviation $N:=\frac14n^{r-1}$
this yields $f_0=4nN$ vertices and $f_1=4rnN$ edges for $\pi(P_n^{2r})$.

The products $(C_n)^r$ have $P:=rn^{r-1}=4rN$ polygon $2$-faces.
In the projection, all these are preserved,
in addition to some of the quadrilateral $2$-faces.

The projected polytope has two types of facets:
There are ``prism'' facets, which involve two of the
polygons, as well as ``cube'' facets,
which in $(C_n)^r$ arise as products of three edges and
$r-3$ vertices, but contain no polygon $2$-faces.
Thus each prism facet is bounded by two polygons,
and each polygon lies in two prism facets.
Hence there are $P=4rN$ prism facets, as well as
some number $C\ge0$ of cube facets.

Now double counting of ridges yields
$6C+(n+2)P=2f_2$. Thus with the Euler equation we get
$C=\frac14(r-2)n^r=(rn-2n)N$.
Finally, counting the vertex-facet incidences according to facets
yields $f_{03}=8C+2nP=(8rn-16n+8rn)N$.
\end{proof}

\begin{corollary}
For each $\eps>0$ there is a $4$-polytope 
whose fatness is larger than $9-\eps$ and
whose complexity is larger than $16-\eps$.
\end{corollary}

\begin{proof}
If we write the flag vector of $\pi(P_n^{2r})$ as
\[
(\tfrac4r, 4, 5-\tfrac3r+\tfrac4n, 1-\tfrac2r+\tfrac4n; 
              16-\tfrac{16}r)\cdot \tfrac14 r\,n^r
\]
then clearly fatness approaches $9$ and complexity 
approaches $16$, for $r,n\rightarrow\infty$.
\end{proof}

\begin{small}

\end{small}
\end{document}